\documentclass{amsart}

\usepackage{latexsym,amsmath,amsfonts,amscd,amssymb, amsthm}
\theoremstyle{plain}
\newtheorem{theorem}{Theorem}[section]

\newtheorem{proposition}{Proposition}[section]
\newtheorem{corollary}{Corollary}[section]

\theoremstyle{definition}
\newtheorem{definition}{Definition}[section]
\newtheorem{example}{Example}[section]

\theoremstyle{remark}
\newtheorem{remark}{Remark}[section]

\title{On the Classification of Rational Sphere Maps}

\author{John P. D'Angelo}

\address{Dept. of Mathematics, Univ. of Illinois, 1409 W. Green St., Urbana IL 61801}

\email{jpda@math.uiuc.edu}

\begin{document}

\begin{abstract} We prove a new classification result for (CR) rational maps from the unit sphere
in some ${\mathbb C}^n$ to the unit sphere in ${\mathbb C}^N$. To so so, we work at the level
of Hermitian forms, and we introduce ancestors and descendants.

{\bf AMS Classification Numbers}: 32H35, 51F25, 32M99.


\end{abstract}

\maketitle

\section{Introduction}

There is considerable literature on proper holomorphic mappings between unit balls 
in possibly different dimensional complex Euclidean spaces.
See for example [D1], [D2], [DHX],  [Fa1], [Fa2], [Fo], [H2], [HJ], [JZ], [L], [LP]  and their references.  
By a well-known result of Forstneri$\check{c}$ [Fo], when the domain dimension
is at least $2$, and a proper map $f$ between balls is assumed sufficiently differentiable at the boundary sphere, then $f$ is a rational function. By a result of Cima-Suffridge [CS], $f$ has no singularities on the sphere. Thus $f$ maps the unit sphere in the domain to the unit sphere in the target.
We write ${\mathcal R}(n,N)$ for the collection of rational maps sending the unit sphere in ${\mathbb C}^n$ to the unit sphere in ${\mathbb C}^N$; we allow domain dimension $1$ and we include constant maps. We write ${\mathcal R}^*(n,N)$
for the non-constant maps in ${\mathcal R}(n,N)$.
Despite many papers on this topic, the collection
of rational sphere maps is not well understood when $N$ is large relative to $n$. 

In this paper we introduce two new ideas. 
First, we define ancestors and descendants of the Hermitian forms corresponding to rational sphere maps. 
We show that every rational sphere map is an {\it ancestor} of a {\it final descendant}.
We then fix the denominator and the degree of the numerator, and provide a partial classification
of the final descendants. When the rational sphere map is a polynomial, we recover (in different language)
a result of the author. See Theorem 4.1. When the denominator is of first degree, we give a decisive result classifying the
possible numerators of final descendants.  See Theorem 5.2. 
The general situation (Theorems 5.1 and 7.1) uses similar ideas but it is harder to state the results in simple language. 
In all cases there is a canonical subspace of the target space of a final descendant, on which we give complete information.
In the polynomial case this space is the full target space. When the denominator is of degree
$1$, this subspace, although proper, tells the full story. Additional complications arise when the degree of $q$ exceeds $1$.

We also study an invariance property of the Hermitian form associated with the final
descendant of a rational sphere map, proving the following result.
If $f$ is the final descendant of a rational sphere map, and its associated Hermitian
form is invariant under the circle action $z \mapsto e^{i\theta}z$, then either $f$
is a polynomial (and hence $Uz^{\otimes m}$ for $U$ unitary), or its numerator and denominator
have the same degree, thereby simplifying the classification. 
In [DX1] and [DX2], the author and Xiao have systematically studied analogues of this result
when the Hermitian form is invariant under general subgroups of the automorphism group of the unit ball.

Section 2 summarizes known results on the complexity of rational sphere maps. Sections 3 and 4
show how Hermitian forms arise in this setting; we introduce ancestors and descendants in Section 4. 
Section 5 describes in detail the relevant linear system of equations involving the inner products
of vector-valued homogeneous polynomials. We prove Theorem 5.2 there.
We prove the invariance result in Section 6 
and the results about higher order denominators in Section 7.

The author acknowledges support from NSF Grant DMS 13-61001. He thanks Ming Xiao
and Jiri Lebl for useful discussions about related ideas and the referee for several
valuable comments.

\section{A summary of known results}

Let ${\mathbb C}^n$ denote complex Euclidean space with norm $|| \ ||$ and inner product $\langle \ , \ \rangle$.
We recall that the holomorphic automorphism group of the unit ball $B_n$ in ${\mathbb C}^n$ consists
of linear fractional transformations of the form $U \phi_a$, where $||a||<1$, where $U$ is unitary, and 
$$\phi_a(z) = {a - L_az \over 1 - \langle z,a \rangle} $$
for a linear map $L_a$ depending on $a$. With $s^2 = 1 - ||a||^2$, one has
$$ L_a(z) = {\langle z, a\rangle a \over s+1} + s z. \eqno (1) $$

Constants and ball automorphisms provide the simplest examples
of {\it rational sphere maps}. 
We consider rational functions $f:{\mathbb C}^n \to {\mathbb C}^N$ such that the image of the unit sphere in the domain
lies in the unit sphere in the target. Thus $||f(z)||^2=1$ whenever $||z||^2=1$.
We write ${\mathcal R}(n,N)$ for the collection of such maps.
These maps are functions of $z=(z_1,...,z_n)$ and independent of the ${\overline z}$ variables.
We assume that the rational map is reduced to lowest terms and (when $n\ge 2$) 
that the constant term in the denominator equals $1$. We write ${\mathcal R}^*(n,N)$ for the non-constant maps
in ${\mathcal R}(n,N)$.

We use the following notations throughout this paper. We write $f\oplus 0$ for the map to ${\mathbb C}^N$ given by
$z \to (f(z), 0)$. More generally $f \oplus g$ denotes the orthogonal sum of $f$ and $g$. Next, suppose
$f:{\mathbb C}^n \to {\mathbb C}^N$ is a polynomial of degree $d$. We write
$$ f = \sum_{k=\nu}^d f_k $$
to denote its decomposition into vector-valued homogeneous polynomials. See also Section 3.

\begin{definition} Assume $f,g \in {\mathcal R}(n,N)$. They are {\bf spherically equivalent} if there are automorphisms $\phi$ of the domain ball and $\chi$ of the target ball such that $f = \chi \circ g \circ \phi$. \end{definition}

\begin{definition} Assume $f,g \in {\mathcal R}^*(n,N)$. They are {\bf homotopically equivalent in target
dimension} $N$ if there is a one-parameter family $H_t$ such that
\begin{itemize}
\item $H_t \in {\mathcal R}^*(n,N)$ for each $t \in [0,1]$.
\item $H_0 = f\oplus 0 $ and $H_1 = g \oplus 0$.
\item The coefficients of $H_t$ depend continuously on $t$. 
\end{itemize}
\end{definition}

Spherical equivalence implies homotopy equivalence, but the converse is false.

The second item in Definition 2.2 may require some comment. Consider the family of maps
$ z \to \left(\cos(\theta)z, \sin(\theta)z^2\right)$ from ${\mathbb C}$ to ${\mathbb C}^2$.
This family shows that $(z,0)$ and $(0,z^2)$ are homotopic in target dimension $2$; since the unitary group
is path-connected, it follows that $(z,0)$ and $(z^2,0)$ are homotopic in target dimension $2$.
But $z$ and $z^2$ are {\bf not} homotopic in target dimension $1$. It can be shown that
any $f,g \in {\mathcal R}(n,N)$ are homotopic in target dimension $N+1$ if we allow constant maps in the homotopy,
and in target dimension $N+n$ if we allow only non-constant maps.

The following theorem (See [D2]) shows that ${\mathcal R}(n,N)$ is a large set when $N$ is large.
It also suggests that classification by target dimension might be impossible.

\begin{theorem} Let ${p \over q}:{\mathbb C}^n \to {\mathbb C}^K$ be an arbitrary rational map, reduced to lowest terms,
 with $||{p \over q}(z)||^2 < 1$ when $||z||^2 \le 1$. 
Then there is a polynomial $h:{\mathbb C}^n \to {\mathbb C}^M$ such that $ {p \oplus h \over q} \in {\mathcal R}(n,M+K)$.
\end{theorem}

In Theorem 2.1, no bounds on $N=M+K$ or deg$(h)$ in terms of $K,n$ and deg$(p)$ alone are possible.
Theorem 2.1 implies that every polynomial which is non-zero on the unit sphere is the denominator for a rational sphere map,
reduced to lowest terms. Thus classification of rational sphere maps requires regarding something as fixed.
Up to now most authors have fixed the target dimension. In this paper we will proceed differently by fixing the denominator.

We summarize well-known results in the next several propositions and theorems.

\begin{proposition} The set ${\mathcal R}(n,N)$ consists only of constants when $N<n$. \end{proposition}
\begin{proof} Since $n \ge 2$ here, a non-constant rational sphere map extends to a proper map $f$ of the unit balls.
The inverse image of a point would be a compact positive dimensional complex analytic subvariety of the ball,
which is impossible. Thus in this case each rational sphere map is a constant. \end{proof}

The next result is a special case of much more general results giving circumstances when proper
holomorphic self-maps are necessarily automorphisms. See [P] and the recent work [J].

\begin{proposition} For $n \ge 2$, the set ${\mathcal R}(n,n)$ consists only of constants and 
linear fractional transformations which are automorphisms
of the unit ball. \end{proposition}

\begin{proposition} The set ${\mathcal R}(1,1)$ is the set of functions given by
$$ e^{i\theta} z^{-m} \prod_{j=1}^K { z - a_j \over 1- {\overline a_j} z}, \eqno (2) $$
where $|a_j| \ne 1$ and $m,K \ge 0$. \end{proposition}
Note that a proper holomorphic map of the unit disk
is a finite Blaschke product; each $|a_j| < 1$ and $m=0$. Factors in (2) where $|a_j|=1$ are constant, and hence
omitted. Only in one dimension
are there non-constant rational sphere maps that are not proper holomorphic maps of the ball, corresponding
to the term $z^{-m}$ or to factors with $|a_j|>1$ in (2).

One of the author's aims has been to view these three Propositions as part of a unified theory.
Before getting to those ideas, we continue with our summary.

In the next two results the authors proved stronger theorems than we state here, as they considered proper
maps between balls with some regularity at the boundary.

\begin{theorem}[Faran] Consider rational sphere maps ${\mathcal R}(n,N)$.
\begin{itemize} 
\item Each $f \in {\mathcal R}(2,3)$ has degree at most $3$. There are four spherical equivalence classes
of proper maps in ${\mathcal R}(2,3)$. See [Fa1].
\item Assume $2 \le n \le N \le 2n-2$. Each $f \in {\mathcal R}(n,N)$ has degree at most $1$. There is
one spherical equivalence class of proper maps in ${\mathcal R}(n,N)$. See [Fa2].
\end{itemize} \end{theorem}

\begin{theorem}[Huang-Ji] For $n \ge 3$, each $f \in {\mathcal R}(n,2n-1)$ has degree at most $2$. 
There are two spherical equivalence classes of proper maps in ${\mathcal R}(n,2n-1)$. See [HJ]
and also [H1],[H2] for related results.  \end{theorem}

\begin{theorem}[Lebl] In source dimension at least $2$,
each quadratic rational sphere map is spherically equivalent to a quadratic monomial map. See [L].
\end{theorem}

We remark that Theorem 2.4 fails in one dimension; the simplest example is given by ${1 \over z^2}$.

\begin{theorem} For $N\ge 2n$ and all $n$, the maps in ${\mathcal R}^*(n,N)$ lie in
infinitely many spherical equivalence classes. In fact there is a one-parameter family $H_t$
of quadratic polynomials where each $H_t$ lies in a different spherical equivalence class. \end{theorem}

See [DL1] for the following stronger statement and finiteness theorem. 

\begin{theorem}[D'Angelo-Lebl] Let $H_t$ be a homotopy of rational proper maps between balls.
Then either all the maps are spherically equivalent or there are uncountably many
spherical equivalence classes. \end{theorem}

\begin{theorem}[D'Angelo-Lebl] For $n\ge 2$ and each $N$, the maps in ${\mathcal R}^*(n,N)$ lie in
finitely many homotopy equivalence classes.\end{theorem}

We close this section by discussing the degree estimate conjecture made by the author many years ago.

\medskip

\noindent {\bf Conjecture}. Assume $f \in {\mathcal R}(n,N)$. The following sharp bounds hold:
\begin{itemize}
\item If $n =2$, then $deg(F) \le 2N-3$.
\item If $n\ge 3$, then $deg(F) \le {N-1 \over n-1}$. 
\end{itemize}

By Proposition 2.3, there is no bound on the degree of elements in ${\mathcal R}(1,1)$ and hence
in ${\mathcal R}(1,N)$. The conjecture is known for monomial maps; see [DKR] for $n=2$ and  [LP] for $n \ge 3$.
The proofs are rather complicated in both cases. There are explicit
examples with the given degrees, and hence the conjecture would be sharp
if proved. The following non sharp-bound is known. See [DL2].

\begin{theorem}[D'Angelo-Lebl]  If $n\ge 2$, and $f \in {\mathcal R}(n,N)$ is of degree $d$, then 
$$ d \le {N(N-1) \over 2(2n-3)}.  \eqno (3) $$ \end{theorem}

The proof of (3) relies on a degree estimate proved by Meylan [M] when $n=2$.

\section{Hermitian forms}

The main tool in this paper is Hermitian forms. Let $W(n,d)$ denote the complex vector space
of polynomials of degree at most $d$ in $n$ variables, and
let $V(n,d)$ denote the subspace of homogeneous polynomials of degree $d$.
(Of course we also include the zero polynomial.)
It is possible to homogenize by adding a variable and then to work with $V(n+1,d)$ instead of $W(n,d)$,
but it makes no essential difference. 

We make $V(n,m)$ into an inner product space by decreeing that the distinct monomials
are orthogonal and that $||z^\alpha||_V ^2 = {m \choose \alpha}$, the multinomial coefficient.
We will often work with ${\mathbb C}^N$-valued homogeneous polynomials. We use the following abbreviated notation.
Suppose, with coefficients $c_\alpha \in {\mathbb C}^N$, we have
$$ p(z) = \sum_{|\alpha|=m} c_\alpha z^{\alpha}. $$
We write $p(z) = L(z^{\otimes m})$ where $L:V(n,m) \to {\mathbb C}^N$ is a linear map. Thus $z^{\otimes m}$
amounts to a list of the monomials forming an orthonormal basis.

\begin{example} Suppose $p(z_1,z_2) = (z_1^2, \sqrt{2}z_1z_2, z_2^2)$. Then
$$  p(z_1,z_2) = \begin{pmatrix}  1 & 0 & 0 \cr 0 & 1 & 0 \cr 0 & 0 & 1 \end{pmatrix} \ \begin{pmatrix} z_1^2 \cr \sqrt{2} z_1z_2 \cr z_2^2 \end{pmatrix} = L(z^{\otimes 2}).                      $$\end{example}

Next we discuss positivity conditions. In coordinates, a Hermitian form on $W(n,d)$ can be written
$$ r(z,{\overline z}) = \sum_{|\alpha|,|\beta| \le d}  c_{\alpha \beta} z^ {\alpha} {\overline z}^\beta $$
where the matrix $(c_{\alpha \beta})$ is Hermitian symmetric.
We make a well-known but crucial comment:
the conditions that $r(z,{\overline z})$ be non-negative as a function 
of $z$ and be non-negative as a Hermitian form differ.
If the form is non-negative definite, then the function is non-negative. If the function is
non-negative, then the form can have some negative eigenvalues.

\begin{example} Put $r(z,{\overline z}) = (|z_1|^2 - |z_2|^2)^2$.
As a function, $r$ is non-negative. The underlying Hermitian form on $V(2,2)$ is diagonal with eigenvalues
$1,-2,1$. \end{example}

We write 
$$ r(z,{\overline z}) = \sum_{|\alpha|, |\beta| \le d}  c_{\alpha \beta} z^ {\alpha} {\overline z}^\beta \succeq 0 \eqno (4) $$
when the matrix $(c_{\alpha \beta})$ has only non-negative eigenvalues, and we use the symbol $\succ$
when all the eigenvalues are positive. We note that (4) holds if and only if $r$
is a Hermitian squared norm; that is, there are holomorphic polynomials $f_j$ of degree at most $d$
such that
$$ r(z,{\overline z}) = \sum_{j=1}^K |f_j(z)|^2 = ||f(z)||^2. $$

Rational sphere maps ${p \over q}$  will correspond to certain Hermitian forms 
$||p||^2 - |q|^2$ with exactly {\bf one}  negative eigenvalue. See (7) below.

The following theorem (See [Q], [CD], [D2]) plays a key role in the proof of Theorem 2.1.

\begin{theorem}[Quillen, Catlin-D'Angelo] Suppose 
$$ r(z,{\overline z}) = \sum_{|\alpha|=|\beta|=m} c_{\alpha \beta} z^ {\alpha} {\overline z}^\beta > 0$$
on the unit sphere. Then there is an integer $d$ such that
$$ ||z||^{2d} r(z,{\overline z}) \succ 0.$$ \end{theorem}

\begin{example} Consider the monomial $\alpha z w$ in two variables. This function is a component
of a polynomial map to some sphere if $|\alpha| < 2$.
To find the target dimension and degree $d+2$ that work we need
an inequality on Hermitian forms:
$$ |\alpha|^2 |z|^2 |w|^2 ( |z|^2 + |w|^2)^d \preceq (|z|^2 + |w|^2)^{d+2}. \eqno (5) $$
For each $k$ we therefore require
$$ |\alpha|^2 {d \choose k} \le {d+2 \choose k+1} $$
and thus for $0 \le k \le d$ we obtain 

$$ |\alpha|^2 \le {(d+2)(d+1) \over (k+1)(d-k+1)} . $$
Assuming $d$ is even, the critical value is when $k = {d \over 2}$.
We get $|\alpha|^2 \le {4(d+1) \over (d+2)}$. Hence (after rewriting), for
$\alpha zw$ to be a component of a map of degree $d+2$, we require 
$$ d \ge {2 |\alpha|^2 - 4 \over 4 - |\alpha|^2}. \eqno (6) $$
By (6), if $|\alpha|$ approaches $2$, then $d$ approaches $\infty$.
A similar situation happens for the target dimension. \end{example}

\section{Hermitian forms and rational sphere maps}

Consider ${p \over q} \in  {\mathcal R}(n,N)$. 
Write $(p;q)$ for the corresponding polynomial map from ${\mathbb C}^n$ to ${\mathbb C}^{N+1}$. 
Thus $q$ is scalar-valued and not $0$ on the closed unit ball.
Without loss of generality we assume $p$ and $q$ have no common factors and that $q(0)=1$.

Given an arbitrary polynomial map $(p;q)$ to ${\mathbb C}^{N+1}$ with no common factors, 
we ask when it corresponds to an element of ${\mathcal R}(n,N)$ or ${\mathcal R}^*(n,N)$.
We want $||p||^2 - |q|^2 = 0$ on the sphere. Equivalently, there are polynomial maps $f,g$ such that

$$ ||p||^2 - |q|^2 = \left(||f||^2 - ||g||^2 \right) \ \left(||z||^2 -1 \right). \eqno (7) $$
We call $||f||^2 - ||g||^2$ the {\it quotient form}. The quotient form is not $0$ when ${p \over q}$
is not a constant. Then formula (7) is equivalent to saying that

$$ ||(f \otimes z) \oplus g ||^2 - ||(g \otimes z ) \oplus f||^2 $$
has signature pair $(N,1)$. Thus,  as a Hermitian form, there are $N$ positive and $1$ negative eigenvalue.
Formula (7) seems easy, but things are quite subtle. 

\begin{definition} Let $g={p \over q} \in {\mathcal R}^*(n,N)$ be a rational sphere map, reduced
to lowest terms and with $q(0)=1$. We define
its associated Hermitian form ${\mathcal H}(g)$ by
$$ {\mathcal H}(g) = ||p||^2 - |q|^2. $$
When $g$ is of degree $d$, we regard ${\mathcal H}(g)$ as a Hermitian form on the vector space
$W(n,d)$. \end{definition} 

Let $g = {p \over q}$ be a rational sphere map and let $\phi_a$ be an automorphism of the ball in the target space.
Write $G={P \over Q} = \phi_a \circ g$. Then we have
$$ ||P||^2 - |Q|^2 = (1- ||a||^2) (||p||^2 - |q|^2). \eqno (8) $$
Thus ${\mathcal H}(G)$ is a constant multiple of ${\mathcal H}(g)$ and
the quotient form of $G$ is $1-||a||^2$ times the quotient form of $g$.
See [L] for applications such as Theorem 2.4. Theorem 6.1 provides an elegant
result about the invariance of ${\mathcal H}(g)$ under a circle action.

We next compute ${\mathcal H}(g)$ when $g$ is the tensor product
of automorphisms. For $||a|| < 1$ we put $c_a = 1-||a||^2$. We write $\rho = ||z||^2-1$ for the defining
equation of the unit sphere and we put $W_j= W_j(z,{\overline z})  = |1 - \langle z,a_j\rangle|^2$

\begin{proposition} Suppose $g = {p \over q}$ is the tensor product
of $K$ automorphisms $\phi_{a_j}$. We assume each $a_j \ne 0$. Write $c_j$ for $c_{\alpha_j}$.
Then we have the following formula for the Hermitian form corresponding to $g$.

$$ ||p||^2 - |q|^2 = \prod_{j=1}^K (c_j \rho + W_j)- \prod_{j=1}^K W_j. \eqno (9) $$
\end{proposition}

\begin{proof} By (8), applied when $g$ is the identity map, the squared 
norm of the numerator of $\phi_{a_j}$ can be written:
$$ c_{j} \rho + W_j.  $$
Note that the squared norm of a tensor product is the product of the squared norms of the factors.
Hence the numerator of ${p \over q}$ is the tensor product of the corresponding numerators and the denominator is the product of the corresponding denominators.
As claimed, we obtain
$$ ||p||^2 - |q|^2 = \prod_{j=1}^K (c_j \rho + W_j)- \prod_{j=1}^K W_j.  $$\end{proof}

Formula (9) defines a polynomial $\sum_{j=1}^K B_j \rho^j$ in the defining function $\rho$. The
coefficients $B_j$ are functions, but satisfy simple formulas such as
$$ B_0 = 0 $$ 
$$ B_1 = \sum_j c_j \prod_{k \ne j} W_k$$
$$ B_2 = \sum_{j \ne k} c_j c_k \prod_{l \ne j,k} W_l$$
$$ B_K = \prod_{j=1}^K c_j. $$
These formulas indicate the symmetry of the result in the points $a_j$. 

\begin{definition} Let $r = ||p||^2 - |q|^2$ and $s= ||f||^2 - |q|^2$ be Hermitian forms 
with the same negative term.
We say that $r$ is a {\bf first ancestor} of $s$ or that $s$ is a {\bf first descendant} of $r$ if
$$ s= E(r) = E( ||p||^2 - |q|^2) = ||p||^2 + (||z||^2 -1)\  ||\pi(p)||^2 - |q|^2, \eqno (10) $$ 
and the degree of $s$ equals the degree of $r$.
Here $\pi$ is orthogonal projection onto a nonzero subspace of the target ${\mathbb C}^N$ of ${p \over q}$. For $k \ge 2$,
we say that $s$ is a $k$-th {\bf descendant} of $r$ if it is a first descendant of a $(k-1)$-st descendant
of $r$. We say that $s$ is a { \bf final descendant} of $r$, if we cannot apply the operation $E$
in (10) without increasing the degree. \end{definition}

By (10), the quotient form of $E(r)$ equals the quotient form of $r$ plus $||\pi(p)||^2$.
This fact arises in the proof of a result from [D3], stated below as Theorem 4.2.

We clarify a notational issue. Assume $p=(p_1,..., p_N)$. Then $\pi(p)= (p_{j_1},..., p_{j_k})$ 
can be regarded, after a unitary map, as simply a list of some of the components of $p$. 

Note that $E(r)$ and $r$ are equal on the unit sphere. In particular,
if $r$ corresponds to a rational sphere map, then so does $E(r)$. Furthermore, the denominator is unchanged.
The basic idea of this paper is simple; we start with a rational sphere map ${p \over q}$ and consider its Hermitian
form $r = ||p||^2 - |q|^2$. We apply the operation $E$ until we reach a final descendant.
Then we describe all final descendants. In [D3] this process is called {\bf orthogonal homogenization}.
See Theorem 4.1 below.

The number of positive eigenvalues of the form $r$ is not generally
preserved by the operation in (10). It can increase, decrease, or stay the same.
This situation partially explains why the degree estimate conjecture is difficult.

\begin{example} Consider the Hermitian form (on $W(2,5)$)
$$ |z|^{10} +|w|^{10} + 5 |z|^6 |w|^2 + 5 |z|^2  |w|^4 -1 .$$
It corresponds to the polynomial sphere map $(z^5, \sqrt{5}z^3w, \sqrt{5} z w^2, w^5)$
in ${\mathcal R}(2,4)$. The Hermitian form
$$  |z|^{10} + |w|^{10} + (|z|^2+ |w|^2) (5 |z|^6 |w|^2 + 5 |z|^2  |w|^4 ) -1 $$ 
is a first descendant. The form
$$  (|z|^2+ |w|^2)^5 - 1 = $$
$$ |z|^{10} + |w|^{10} + (|z|^2+ |w|^2) (5 |z|^6 |w|^2) + (|z|^2 + |w|^2)^2 (5 |z|^2  |w|^4 ) -1 $$
is a second (and final) descendant. It corresponds to an element in ${\mathcal R}(2,6)$.
\end{example} 

\begin{remark} The ancestor form in Example 4.1 provides an example of a polynomial sphere mapping invariant
under a cyclic group of order five: $(z,w) \mapsto (\eta z, \eta^2 w)$, where $\eta$ is a $5$-th root of unity.  A map corresponding to the first descendant is not invariant under any non-trivial group.
The map corresponding to the final descendant is invariant under a different representation of 
the cyclic group of order five:  $(z,w) \mapsto (\eta z, \eta w)$. \end{remark}

This new language yields the following reformulation of a result from [D1].

\begin{theorem} Let $p$ be a polynomial sphere map of degree $m$. Then $||p||^2 -1 $ is an ancestor
of $||z||^{2m} - 1 = ||z^{\otimes m}||^2 - 1$.   \end{theorem}

We recall an alternative way to state this result, which focuses on the sphere map rather than on the Hermitian form.

\begin{corollary} Let $p$ be a polynomial sphere map of degree $m$. Then there is a finite
number of tensor product operations $E_1, ..., E_k$ and a unitary map $U$ such that 
$$(E_k \cdots E_1)(p) = U z^{\otimes m}. $$
\end{corollary}

\begin{remark} Let $p$ be a polynomial sphere map. Assume $p$ vanishes to order $\nu$ at $0$ and 
is of degree $m$. Then $z^{\otimes m}$ is a $k$-th descendant of $p$, where $k=m-\nu$.   \end{remark}

We mention a volume inequality from [D3] whose proof uses the operation (10).

\begin{theorem} Let $V_p$ be the volume (with multiplicity counted) of the image of the ball under a polynomial sphere map $p$
of degree $m$. Then

$$ V_p \le {\pi^n m^n \over n!}. \eqno (11) $$ 
Equality occurs in (11) if and only if $p = U z^{\otimes m}$ for some unitary $U$.
\end{theorem}

The proof of Theorem 4.2 relies on the following result. Identify a polynomial map $p$ with the Hermitian
form ${\mathcal H}(p) = ||p||^2-1$. Then the volume of the image of a descendant is greater than the volume of the image of the ancestor.

\section{Equations on inner products}

Consider a ${\mathbb C}^N$-valued rational function ${p \over q}$. We write
$$ p(z) = \sum A_\alpha z^\alpha $$
$$ q(z) = \sum b_\alpha z^\alpha $$
where each $A_\alpha \in {\mathbb C}^N$ and each $b_\alpha \in {\mathbb C}$.
The condition for being a sphere map is a system of linear equations
in the inner products $\langle A_\alpha, A_\beta \rangle$ and the scalars $b_\alpha {\overline b_\beta}$.
For $||z||^2 =1 $ we have: 

$$ ||p(z)||^2 - |q(z)|^2 = \sum_{\alpha, \beta} \left( \langle A_\alpha, A_\beta \rangle - b_\alpha {\overline b_\beta}\right) z^\alpha {\overline z}^\beta = 0.  \eqno (12)  $$
Homogenizing and equating coefficients leads to messy formulas. 

Let $D(n,d) = {d+n-1 \choose n-1}$ denote the dimension of $V(n,d)$.
The following combinatorial result gives the number of linear equations
satisfied by the inner products of
the vector coefficients of a polynomial sphere map of degree $d$ in $n$ variables.
If we regard the denominator of a rational sphere map as known, then the inner products
of the vector coefficients of the numerator satisfy the same number of equations.

\begin{proposition} Let $p(z)= \sum C_\alpha z^\alpha$ denote a polynomial sphere map of degree $d$.
The inner products of
the vector coefficients $C_\alpha$ of $p$ satisfy a linear system of $K=K(n,d)$ equations, where
$$ K(n,d) = D(n,d) \sum_{j=0}^{d-1}D(n,j) \ + \left({D(n,d) (D(n,d)+1) \over 2} \right).  \eqno (13)$$
For fixed $n$, the number $K(n,d)$ is a polynomial of degree $2n-1$ in $d$.
\end{proposition}
\begin{proof} (Sketch) The term $D(n,d)$ in (13) is the dimension of $V(n,d)$.
Each term in the sum is the dimension of $D(n,j)$ for $j <d$. Hence the expression
$D(n,d)$ times the sum results from counting the number of inner products arising when the degree
of homogeneity is $d$ in the $z$ variables and less than $d$ in the conjugated variables.
The other term equals $1+ 2+ \dots + D(n,d)$, which is the number of inner products
arising from terms homogeneous of degree $d$ in both the $z$ variables and in the conjugated variables.
\end{proof}

\begin{example} We have the following results:
\begin{itemize}
\item $K(1,d) = d+1$. 
\item $K(2,d) ={d^3+ 3 d^2 + 4d +2 \over 2}$.
\item $ K(3,d) = {2d^5 + 15 d^4 + 44 d^3 + 69 d^2 + 62 d + 24 \over 24} $.
\end{itemize}
\end{example}

Expanding in terms of homogeneous polynomials is easier. We illustrate with the next example.
There are $34$ equations in the approach where we regard the coefficient vectors
as unknowns, and there are $4$ equations if we regard homogeneous polynomials as unknowns.

\begin{example} Consider the map $F:{\mathbb C}^2 \to {\mathbb C}^{10}$ defined by
$$f(z,w) = A+Bz+Cw+ Dz^2+ Ezw+ F w^2 + Gz^3 + Hz^2 w + I z w^2 + J w^3. $$
We can regard $A,B,C,D,E,F$ are parameters. The inner products involving $G,H,I,J$ are then determined
by the equations in (12) after equating coefficients.

$$ \langle A,G\rangle = \langle A,H \rangle = \langle A, I \rangle = \langle A, J \rangle = 0. $$

$$ \langle Gz^3+Hz^2w+ Izw^2 + Jw^3, Bz+Cw\rangle = $$
$$ - \langle Dz^2+Ezw+Fw^2, A \rangle (|z|^2+|w|^2). $$

$$ \langle Gz^3+Hz^2w+ Izw^2 + Jw^3, Dz^2+Ezw+Fw^2 \rangle = $$
$$ - \langle Dz^2+Ezw+Fw^2, Bz+Cw \rangle (|z|^2+|w|^2) - \langle Bz+Cw, A \rangle (|z|^2+|w|^2)^2. $$
There is one more long equation involving squared norms.
Using the expansion in terms of homogeneous parts the equations become 

$$ \langle p_3, p_0 \rangle = 0 $$

$$ \langle p_3, p_1 \rangle = - \langle p_2,p_0 \rangle (|z|^2 + |w|^2) $$

$$ \langle p_3,p_2\rangle = - \langle p_2,p_1 \rangle (|z|^2+|w|^2) - \langle p_1, p_0 \rangle (|z|^2 + |w|^2)^2  $$

$$ ||p_3||^2 = $$ 
$$ (|z|^2 + |w|^2)^3 - (|z|^2 + |w|^2 ) ||p_2||^2 - (|z|^2 + |w|^2)^2  ||p_1||^2 - (|z|^2 + |w|^2)^3 ||p_0||^2 . $$
We can regard these last four equations as follows. We think of $p_0,p_1,p_2$ as known. All the right-hand sides are then known.
The left-hand sides then tell us the inner products of $p_3$ with each of these lower order terms. 
Expanding each of the homogeneous polynomials in coordinates gives the more complicated system described above.
\end{example}

We approach the analogous equations in the rational case by fixing the denominator and degree of numerator. We will find all maps, ignoring target dimension, and provide a partial classification. Given $q$ with $q(z) \ne 0$ on the sphere, how do we construct 
all possible numerators $p$? Assume the degree of $q$ is $k$. The equations in (12) imply that the degree of $p$
is at least as large as the degree of $q$, and hence we assume $p$ has degree $m+k$ for $m \ge 0$.
The condition $||p||^2 = |q|^2$ on the sphere yields the following  results.

\begin{proposition} Assume $n\ge 2$. Let ${p \over q} \in {\mathcal R}(n,N)$.
Put $q= 1 + \cdots + q_k$. For $||z|| \le 1$ and $0 \le j \le k-1$ we have 
$ |q_{k-j}(z)| < {k \choose j}$.
In particular $ |q_k(z)| < 1$. \end{proposition}
\begin{proof} Choose $z$ with $||z||^2 =1$, and consider the complex line $t \mapsto tz$. The restriction
to this line defines a rational sphere map, with no singularities in the disk,
in one dimension. By homogeneity its denominator is
$$ u(t) = \sum_{j=0}^k q_j(z) t^j.$$
But all the roots $a_j$ of $u$ lie outside the unit disk, and hence we can write
$$ u(t) = \prod_{j=1}^k (1 - {\overline a_j(z)} t), $$
where $|a_j(z)| < 1$ for each $j$. Expanding and estimating gives the result.
\end{proof}

\begin{proposition} Let ${p \over q}$ be a rational sphere map. Assume the degree of $p$ is $m+k$.
Expand $q$ as $q = 1 + \cdots + q_k$ in terms of homogeneous polynomials. Then 
there is a final descendant $||g||^2 - |q|^2$ with
$g_j = 0$ for $j < m$ and 
$$ \langle g_{m+k}, g_m \rangle = q_k ||z||^{2m} = \langle q_k z^{\otimes m}, z^{\otimes m} \rangle. \eqno (14) $$
\end{proposition}
\begin{proof} 
On the sphere we have
$$ \sum_{j,l=0}^{m+k} \langle p_{j},p_{l} \rangle = \sum_{j,l=0}^k q_j {\overline {q_l}}. $$
Replace $z$ by $e^{i\theta}z$ and use homogeneity to get:
$$ \sum_{j,l=0}^{m+k} \langle p_{j},p_{l} \rangle e^{i\theta(j-l)} = \sum_{j,l=0}^k q_j {\overline q_l} e^{i\theta(j-l)}. \eqno (15)$$
We equate Fourier coefficients in this equality of trig polynomials. Put $j-l=b$. Then $-(m+k) \le b \le (m+k)$.
We call $b$ the gap in indices. 

On the sphere, for each $b$ with $-(m+k) \le b \le m+k$,  we get

$$ \sum_{l=0}^{m+k-b} \langle p_{l+b},p_{l} \rangle  = \sum_{l=0}^{k-b} q_{b+l} {\overline q_l}. \eqno (16) $$
Since $q_j = 0$ for $j >k$, we must have $b+l \le k$ in the right-hand side of (16).

The largest gap is when $b=m+k$. Put this value into (16). Both sums have only one term, and we get 
$$ \langle p_{m+k}, p_0 \rangle = q_{m+k}.$$
Since $q$ is of degree $k$, we conclude either that $m=0$ or that $p_{m+k}$ is orthogonal to $p_0$.
If $m=0$, we put $g=p$ and we already have $g_j = 0$ for $j<m$. Furthermore $\langle p_k, p_0 \rangle = q_k$
on the sphere and, by homogeneity, (14) holds. Thus the conclusion holds when $m=0$.

Assume $m > 0$ and that $p_{m+k}$ is orthogonal to $p_0$.
If $p_0 \ne 0$, then we let $V$ be the subspace spanned by $p_0$. Since $p_{m+k}$ is orthogonal to $V$,
we can apply (10) to get a new polynomial $g$ of the same degree $m+k$ with $g_0 = 0$. Thus, whether or not $p_0=0$,
we may assume we have a descendant $||g||^2 - |q|^2$ with $g_0=0$. 

Now the largest gap is when $b=m+k-1$. We apply the same reasoning. If $m=1$, we get 
what we want, since, on the sphere,
$$ \langle p_{m+k}, p_m \rangle =   \langle p_{k+1}, p_1 \rangle = q_k. $$
Homogenizing gives (14). 

If $m>1$, then $p$ vanishes to order at least $2$, and as above, we can apply (10) to fix the degree and
increase the order of vanishing. We can proceed in this way until we get a descendant $||g||^2 - |q|^2$ satisfying
$g = g_m + \cdots + g_{m+k}$. 

Take $b=k$ in (16), with $p$ replaced by $g$. The sum on the right-hand side has only one term (when $l=0$), namely $q_k$.
Hence on the sphere we have
$$ \sum_{l=0}^m \langle g_{l+k}, g_l \rangle = q_k. $$
But now $g_j = 0$ for $j<m$ and the sum on the left-hand side has only one term. 
We conclude that $\langle g_{m+k}, g_m \rangle = q_k$.
Homogenizing gives (14).

We summarize the proof. We expand the numerator and denominator of a rational sphere map
into homogeneous parts. We express the condition for being a sphere map
in terms of inner products. We let the circle act on the sphere by replacing $z$ by $e^{i\theta}z$,
and then equate Fourier coefficients. The resulting identities hold on the sphere. 
We use them in conjunction with the notion of descendant to reduce to the case where
$${g \over q} = {g_m + \cdots + g_{m+k} \over 1 + \cdots + q_k}. \eqno (17) $$
Maps as in (17) satisfy identities such as (18.1-3) below. \end{proof}

Once we have a sphere map ${g \over q}$ satisfying the properties in Proposition 5.1, we can draw several conclusions.
The first conclusion is that there is a canonical non-zero subspace $W$ into which both $g_m$ and $g_{m+k}$ map.
Then (14) provides one of the major parts of Theorem 5.1 below.

We rewrite (14) along with two of the other bihomogenized identities:

$$ \langle g_{m+k}, g_m \rangle = q_k ||z||^{2m} \eqno (18.1) $$

$$ \langle g_{m+k}, g_{m+1} \rangle + \langle g_{m+k-1}, g_m \rangle ||z||^2 = q_k {\overline {q_1}} ||z||^{2m} + q_{k-1} ||z||^{2m+2} \eqno (18.2) $$

$$ \sum_l ||g_{m+l}||^2 \ ||z||^{2k-2l} = \sum |q_l|^2 \ ||z||^{2m+2k-2l}. \eqno (18.3) $$

The condition in (18.1) arises when the gap in the indices is maximal, namely $k$. The condition in (18.2)
arises when this gap is $k-1$, and the condition in (18.3) is when the gap is $0$. We do not write out the intermediate
expressions; they arise from bihomogenizing (16).

These ideas lead to a general result. Recall that $V(n,m)$ is the complex vector space of homogeneous polynomials
of degree $m$ in $n$ variables.

\begin{theorem}  Let $f={p \over q}$ be a rational sphere map. Assume ${\rm deg}(p)= m+k$ and
${\rm deg}(q) = k$.  The following hold:
\begin{itemize}
\item There is a finite number of tensor operations such that
$$ E_s \circ \ ... \circ E_1 (f) = {g_{m} + \dots +  g_{m+k}  \over q }. \eqno (19) $$
In other words $f$ is an ancestor of ${g \over q}$ satisfying (19). 

\item There is a non-zero subspace $W$ of the target space of ${g \over q}$ with $\pi_W(g_m)= g_m$. By the next item,
$W$ is isomorphic to $V(n,m)$.

\item There is an invertible linear map $M:(V,n) \to W $ such that
$$ \pi_W(g_{m+k}) = q_k (M^{-1})^* (z^{\otimes m}).  $$
$$ g_m = M(z^{\otimes m}) \oplus 0. $$

\item There is a complete orthogonal decomposition of the target space described below. 
\end{itemize}
\end{theorem}

Before completing this description, we give a complete analysis when the denominator is of first degree.

\begin{theorem} Let $f={p \over q}$ be a rational sphere map  with linear denominator $1+q_1$ and of 
degree $m+1$. 
The following hold:
\begin{itemize}
\item There is a finite number of tensor operations such that
$$ E_s \circ \ ... \circ E_1 (f) = {g_{m+1} + g_m \over 1+q_1}= {g \over q}. \eqno (20) $$

\item There is a non-zero subspace $W$ of the target space with $\pi_W(g_m)= g_m$. By the next item,
$W$ is isomorphic to $V(n,m)$.

\item There is an invertible linear map $M$ from the space of vector-valued homogeneous
polynomials of degree $m$  to $W$, and a homogeneous mapping $h_{m+1}$, such that 
$$ g_{m+1} = q_1 (M^{-1})^* (z^{\otimes m})  \oplus h_{m+1}. $$
$$ g_m = M(z^{\otimes m}) \oplus 0. $$

\item The Hermitian form defined by
$$ |q_1|^2 \big(||z||^{2m} - ||(M^{-1})^* (z^{\otimes m})||^2\big) + ||z||^2 \big( ||z||^{2m} - ||M(z^{\otimes m})||^2 \big) \eqno (21) $$
is non-negative definite.
 \end{itemize}\end{theorem}

\begin{proof} Proposition 5.1 establishes the first three conclusions of both results. 
We therefore first finish the proof of Theorem 5.2, as it is easier.  We write $q_1 = \langle z, a \rangle$.
$$ {g \over q} = {g_{m+1} + g_m \over 1 + \langle z, a \rangle}. \eqno (22) $$
The expression (22) is a rational sphere map if two conditions are met:

$$ \langle g_{m+1}, g_m \rangle = ||z||^{2m} \langle z,a \rangle \eqno (23) $$

$$ ||g_{m+1}||^2 + ||z||^2 ||g_m||^2 = ||z||^{2m+2} + ||z||^{2m} |\langle z,a \rangle|^2. \eqno (24) $$

If (23) holds, then polarization (regarding $z$ and ${\overline z}$ independently) implies that $\pi_W(g_{m+1})$
is divisible by $q_1$. Since $g_m$ and $g_{m+1}$ are homogeneous, we can then find linear maps $M$ and $L$ such that

$$ g_m = M(z^{\otimes m}) \oplus 0 $$

$$ g_{m+1} = \big( \langle z, a \rangle L(z^{\otimes m}) \big) \oplus  h_{m+1}. $$
Formula (23) forces $L= (M^{-1})^*$.  
Here $h_{m+1}$ is a vector-valued polynomial, the orthogonal projection of $g_{m+1}$ onto the orthogonal complement of $W$.
We need to make (21) hold as well. Solving for the unknown $||h_{m+1}||^2$ gives

$$ ||h_{m+1}||^2 = ||z||^2 \big(||z||^{2m} -  ||M(z^{\otimes m})||^2 \big) + |\langle z, a\rangle|^2 \big(||z||^{2m} - ||L(z^{\otimes m})||^2 \big). \eqno (25) $$
Since Hermitian squared norms
correspond to non-negative definite Hermitian forms,
the form on the right-hand side of (25) must be non-negative definite. \end{proof}

We have proved Theorem 5.2 and the first three parts of Theorem 5.1.

\begin{corollary} Let ${p \over q}$ be a rational sphere map of degree $m+1$ and assume
the degree of $q$ is $1$. Then the final descendant of $||p||^2 - |q|^2$ is completely
determined by the linear map $M$ satisfying (21). Maps $M$ and $M'$ give the same descendant
if and only if there is a unitary map $U$ on $W$ such that $M'=UM$. \end{corollary}

\begin{proof} Let $W$ be the subspace in the Theorem. The final descendant is given by
$$  ||(M+q_1 (M^{-1})^*) (z^{\otimes m})||^2 + ||h_{m+1}||^2 - |q|^2.  \eqno (26)$$
The second term in (26) is determined by  (21). The first term is determined by $q_1$ and $M$.
We note that the first term is unchanged if we replace $M$ by $UM$, where $U$ is unitary on $W$, because 
$$ UM + q_1 ((UM)^{-1})^* = UM + q_1 (M^{-1}U^{-1})^* $$
$$ = UM + q_1 U (M^{-1})^* = U \left(M+ q_1 (M^{-1})^* \right) \eqno (27) $$
We have used $U = (U^{-1})^*$. Taking norms of (24) on $W$ gives the result.
\end{proof}

The final descendant of a polynomial sphere map is completely determined by its degree.
The final descendant of a rational sphere map with denominator of first degree
is determined by the subspace $W$ and the isomorphism $M:V(n,m) \to W$.

We wish to interpret (21), or equivalently (25),  in terms of the spectrum of $M$.
Recall that distinct monomials are orthogonal and that $||z^\alpha||^2 = {m \choose \alpha}$.
After identifying $V(n,m)$ with the subspace $W$, we can express our inequalities in terms of eigenvalues of $M$.
Assume $Mv = \lambda v$ and hence $||L(v)||^2 = {||v||^2 \over |\lambda|^2} $.

After some computation with (25) we obtain
$$ ||h||^2 = ||z||^{2m} (1-|\lambda|^2) \left( ||z||^2 - {|\langle z, a\rangle|^2 \over |\lambda|^2}\right). \eqno (28) $$
Using the Cauchy-Schwarz inequality, the condition becomes $||a|| \le |\lambda| \le 1$.

\begin{remark} In this situation, but not for general denominators, a simplification occurs.
On the right-hand side of (28), we multiply a form on $V(n,1)$ by $||z||^{2m}$. But, unlike in the situation of
Theorem 3.1, this term doesn't impact positive semi-definiteness. A form on $V(n,1)$ is a Hermitian squared
norm if and only if it is non-negative as a function. 
\end{remark}

The inequality in (21) is fundamental. It constrains the map $M$; we note that as $a$ tends to $0$,
$M$ is less constrained, as it may have eigenvalues of even smaller modulus.
This seems at first a bit counter-intuitive; when $a=0$
we get a polynomial map. The polynomial we get, however, is not a final descendant.
We illustrate this issue in the one-dimensional case; similar examples apply in general.

\begin{example} For $|a|<1$ and $0 \le \theta \le {\pi \over 2}$, consider the family of sphere maps
$$ \left(\cos(\theta) z {a-z \over 1- {\overline a}z}, \sin(\theta)z \right) =  {1 \over 1- {\overline a}z} \left(\cos(\theta) (az-z^2) , \sin(\theta)(z - {\overline a} z^2) \right).$$
For $0 < \theta < {\pi \over 2}$, these maps are in the form ${p_2 + p_1 \over 1+q_1}$. 
A simple computation verifies that $\langle p_2, p_1 \rangle = q_1 |z|^2 = -{\overline a}|z|^2 z$. 
As $a$ tends to $0$, the limit map is
$$ (-\cos(\theta) z^2, \sin(\theta)z). $$
We see that the coefficient $\mu$ of $z$ (the analogue of $M$) lies anywhere in $0 < |\mu| < 1$. 
\end{example}

In Theorem 4.1, when $p$ is a polynomial, and hence $q=1$, the final Hermitian form is completely determined. 
In Theorem 5.1, when $q$ is of first degree, the final Hermitian form is not completely determined.
It is determined by the isomorphism $M$, whose spectrum lies in the annular region
$$ ||a|| \le |\lambda| \le 1. $$
The target space for ${p \over q}$
contains an isomorphic copy of $V(n,m)$. The orthogonal complement contains {\it excess};
it is determined up to a unitary. 

\begin{remark} Theorem 5.2 provides information
for sphere maps of degree $1$, and even about formula (1) for automorphisms. 
A rational sphere map of degree $1$ 
can be written ${p_0 + p_1 \over 1 + q_1}$. The bihomogenized version of the identities in (18) become
$$ ||p_0||^2 \ ||z||^2 + ||p_1||^2 = ||z||^2 + |q_1|^2, \eqno (29.1 )$$
$$ \langle p_1, p_0 \rangle = q_1 = \langle z, a \rangle. \eqno (29.2) $$
When $m=0$, equation (22) yields
$$ ||h(z)||^2 = ||z||^2(1- ||p_0||^2) + |\langle z,a \rangle|^2 \left( 1- {1 \over ||p_0||^2}\right). \eqno (30) $$
Thus the right-hand side of (30) must be non-negative. 
There is no excess if and only if $h(z)=0$. This condition 
occurs when $||p_0||=1$, and the map is constant. The opposite end of
the scale is when $||a||= ||p_0||$. Up to a unitary we may put $p_0=a$.
Then $p_1 = \langle z,a \rangle {a \over ||a||^2} \oplus h$. 
In this case we write $p_1(z) = L_a(z)$ and discover that $\langle L_a(z), a \rangle = \langle z,a\rangle$.
Therefore $L_a(z) = v \langle z,a \rangle + c z$ for a constant vector $v$ and constant scalar $c$.
Plugging this ansatz in (29.1) and (29.2) determines $v$ and $c$, thereby yielding formula (1).
 \end{remark}

\section{Invariance under a circle action}

Let $g$ be a rational sphere map, and let $f={p \over q}$ be its final descendant. 
We are interested in invariance of the form ${\mathcal H}(f)$ under subgroups of the automorphism group
of the ball. In this section  we prove one such result.

\begin{theorem} Let $f={g \over q}$ be the final descendant of a rational sphere map of degree $d$.
Suppose that ${\mathcal H}(f)$ is invariant under the map $z \mapsto e^{i\theta}z$.
Then either
\begin{itemize}
\item $f$ is a polynomial and hence $f=U z^{\otimes d}$, or
\item Both $g$ and $q$ are of degree exactly $d$.
\end{itemize}

\end{theorem}
\begin{proof} Expand $g$ and $q$ in terms of homogeneous polynomials.
We may assume that (17) holds, where $d=m+k$. Note that (18.1) then holds as well.
We write
$$ {\mathcal H}(f) = ||g||^2 - |q|^2 = ||\sum g_j||^2 - |\sum q_j|^2. \eqno (31) $$
Expanding (31) yields many terms. When we replace $z$ by $e^{i\theta}z$ in this expansion,
there is only one term of the form $c_k e^{ik\theta}$, namely 
$$ c_k = \langle g_{m+k}, g_m \rangle - q_k. $$
By (18.1), however, we can write this term as
$$ (||z||^{2m}-1) q_k. \eqno (32) $$
If ${\mathcal H}(f)$ is invariant, then this term must vanish. It follows that either
$m=0$, in which case $g$ and $q$ are both of degree $k$, or that $q_k=0$. If $q_k = 0$ then
(18.1) implies that $g_m$ is orthogonal to $g_{m+k}$. Suppose $g_m$ is not $0$.
In this case, however, we can apply the tensor
product operation on the space spanned by $g_m$, and obtain a sphere map $E(f)$ which is still
of degree $m+k$. Since $f$ is assumed to be a final descendant, we get a contradiction.
Hence $g_m = 0$. We are now in the same situation as before, except that we have increased
the order of vanishing of the numerator and lowered the degree of the denominator. We can
proceed in this fashion to establish that $g_j =0$ for $m \le j < m+k=d$ and that $q$ is of degree $0$.
We conclude that $g=g_{m+k}$ and $q=1$. Therefore $f$ is a homogeneous polynomial sphere mapping 
of degree $d$. Corollary 4.1 implies the conclusion in this case.

We repeat the idea; invariance under the circle action forces a certain term to vanish.
That term is $(||z||^{2m}-1) q_k$. When the first factor vanishes, we get a map
whose numerator and denominator have the same degree. When the second factor vanishes, we have lowered
the degree of $q$ and increased the order of vanishing of $g$. Because $f$ is a final descendant, 
invariance allows us to repeat the process until we obtain $q=1$.  Thus $f$ is a polynomial
and also a final descendant; Corollary 4.1 implies the desired conclusion.
\end{proof}

\section{Denominators of higher degree}

The part of Theorem 5.1 already proved shows that there is a natural subspace $W$ of the target space
isomorphic to $V(n,m)$. Assume $n\ge 2$. Here $m$ is the order of vanishing of the (final descendant) map ${g \over q}$ at $0$, and hence also the order of vanishing of $g$ there. If ${g \over q}$ is a polynomial, and hence the final descendant
form is $||g||^2 -1 $, then $W$ is the full target space.
If the degree of $q$ is $1$, then $W$ is a proper subspace
of the target corresponding to the final descendant Hermitian form $||g||^2 - |q|^2$. But the map is determined up to a unitary
map by the isomorphism $M$. When $q$ has degree at least $2$, additional subspaces arise because of (18.1-3).

We analyze this fully in degree two.
Assume $q$ has degree $2$ and $p$ has degree $m+2$.
Theorem 5.1 shows that
the final descendant of ${p \over q}$ will have the form
$$ ||g_m + g_{m+1} + g_{m+2}||^2 - |q|^2. $$
Also $$ \langle g_{m+2}, g_m \rangle = q_2 ||z||^2. $$
Again there is a subspace $W$ and a linear map $M:V(n,m) \to W$ such that
$\pi_W(g_m) = g_m$ and 
$ g_m = M(z^{\otimes m})$ and $\pi_W(g_{m+2}) = q_2 (M^{-1})^* (z^{\otimes m})$.
We recall the additional equations involving $q_1$ and $g_{m+1}$.
We have
$$ \langle g_{m+2}, g_m \rangle = q_k ||z||^{2m} \eqno (33) $$

$$ \langle g_{m+2}, g_{m+1} \rangle + \langle g_{m+1}, g_m \rangle ||z||^2 = q_k {\overline {q_1}} ||z||^{2m} + q_{k-1} ||z||^{2m+2} \eqno (34) $$

$$ ||g_{m+2}||^2 + ||g_{m+1}||^2 \ ||z||^2  + ||g_m||^2 \ ||z||^4 = |q_2|^2 \ ||z||^{2m} + |q_1|^2 \ ||z||^{2m+2} + ||z||^{2m+4}. \eqno (35) $$

We return to our study of ${\mathcal R}(n,N)$ with a  general denominator.
Assume that $||p||^2 - |q|^2$ is a final descendant. Then 
$$ \sum ||p_{m+j}||^2 \ ||z||^{2k-2j} = \sum |q_j|^2 \ ||z||^{2m+2k-2j} \eqno (36) $$
holds by (18.3). Hence there is (a large dimensional!) unitary matrix $U$ with

$$ \begin{pmatrix} p_{m+k} \cr ... \cr p_{m+j} \otimes z^{\otimes(k-j)} \cr ... \cr p_m \otimes z^{\otimes k}\end{pmatrix}  = 
U \ \begin{pmatrix} q_{k}\ z^{\otimes m} \cr ... \cr q_{j} \ z^{\otimes(m+ k-j)} \cr ... \cr q_0 \ z^{\otimes (m+k)} \end{pmatrix}. $$
Thus we could put $p_{m+j} = q_j z^{\otimes m} \otimes v_j  \oplus f_{m+j}$.

For simplicity we choose $U$ in a simple way.
The method does not construct all possible numerators. One must take additional subspaces into account.

\begin{theorem} Let $f$ be a rational sphere map of degree $m+k$ with denominator of degree $k$.
Then there is a finite number of tensor operations such that
$$ E_s \circ \ ... \circ E_1 (f) = {\sum_{j=0}^k p_{m+j} \over \sum_{j=0}^k q_j}. \eqno (37) $$
Given $q$ we can construct $p$ as follows:
Choose $v_j$ with $\langle v_i, v_j \rangle = 1$ for $i \ne j$ 
such that the Hermitian form defined by
$$ ||z||^{2m} \big( \sum_j |q_j|^2 (1-||v_j||^2) ||z||^{2k-2j}\big)  \eqno (38) $$
is positive semi-definite. There are vector-valued maps $f_{m+j}$ such that 
$$ p_{m+j} = q_j (z^{\otimes m} \otimes v_j) \oplus f_{m+j}.$$
\end{theorem}

\begin{remark} Observe that, although the sum in (38) is non-negative as a function,
the terms in the sum in (38) can be of both signs. Hence, without the factor $||z||^{2m}$,
the resulting Hermitian form need not be positive semi-definite.
Hence this factor is typically needed as in Theorem 3.1. \end{remark}

\medskip

\section{bibliography}

\medskip

[CD] D. W. Catlin and J. P. D'Angelo,
A stabilization theorem for Hermitian forms and 
applications to holomorphic mappings, {\it Math Research Letters} 3 (1996), 149-166. 

\medskip

[CS] J. A. Cima and T. J. Suffridge, Boundary behavior of rational proper maps. {\it Duke Math. J.} 60 (1990), no. 1, 135¨C138.

\medskip

[D1] J. P. D'Angelo,  Several Complex Variables and the Geometry of Real Hypersurfaces,
CRC Press, Boca Raton, Fla., 1992.

\medskip
[D2] J. D'Angelo, Proper holomorphic mappings,
positivity conditions, and isometric imbedding, {\it J. Korean Math Society}, May 2003, 1-30.

\medskip

[D3] 
J. P. D'Angelo, Hermitian analysis. From Fourier series to Cauchy-Riemann geometry, Cornerstones, Birkh\"auser/Springer, New York, 2013. 

\medskip

[DL1] J. P. D'Angelo and J. Lebl, Homotopy equivalence for proper holomorphic mappings, {\it Adv. Math.} 286 (2016), 160-180.

\medskip

[DL2] J. P. D'Angelo and J. Lebl, On the complexity of proper mappings between balls, {\it Complex Variables and Elliptic Equations},
Volume 54, Issue 3, Holomorphic Mappings (2009), 187-204.

\medskip

[DHX] J. D'Angelo, Z. Huo, and M. Xiao, Proper holomorphic maps from the unit disk to some unit ball, 
{\it Proc. A. M. S.}, Vol. 145, No. 6, (2017), 2649-2660.

\medskip

[DX1] J. P. D'Angelo and M. Xiao, Symmetries in CR complexity theory, {\it Advances in Math} 313 (2017), 590-627.

\medskip

[DX2] J. P. D'Angelo and M. Xiao, Symmetries and regularity for holomorphic maps between balls, (submitted).

\medskip

[DKR] J. P. D'Angelo, S. Kos, and E. Riehl,  A sharp bound for the degree of proper monomial mappings between balls, 
{\it J. Geom. Anal.} 13 (2003), no. 4, 581-593.

\medskip

[Fa1] J. J. Faran, The linearity of proper holomorphic maps between balls in the low codimension case,
 {\it J. Diff. Geom.} 24 (1986), no. 1, 15-17. 

\medskip

[Fa2] J. J. Faran, Maps from the two-ball to the three-ball, {\it Invent. Math.} 68 (1982), no. 3, 441-475.

\medskip

[Fo] F. Forstneri$\check{c}$,  Extending proper holomorphic maps of positive codimension,
{\it Inventiones Math.}, 95 (1989), 31-62.

\medskip

[H1] X. Huang, On a linearity problem for proper maps between balls in complex spaces
of different dimensions, {\it J. Diff. Geometry} 51 (1999), no 1, 13-33.

\medskip

[H2] X. Huang, On a semi-rigidity property for holomorphic maps, 
{\it Asian J. Math.} 7 (2003), no. 4, 463–492.

\medskip

[HJ] X. Huang, X., and S. Ji, 
Mapping ${\bf B}_n$ into ${\bf B}_{2n-1}$, {\it Invent. Math.} 145 (2001), 219-250.
13-36.

\medskip 
[J] J. Janardhanan,
Proper holomorphic mappings of balanced domains in $C^n$, 
{\it Math. Z.} 280 (2015), no. 1-2, 257-268.

\medskip
[JZ] S. Ji and Y. Zhang, Classification of rational holomorphic maps from $B^2$ into $B^N$ with degree $2$,
{\it Sci. China Ser. A } 52 (2009), 2647-2667.

\medskip

[L] J. Lebl,  Normal forms, Hermitian operators, and CR maps of spheres and hyperquadrics,
{\it Michigan Math. J.} 60 (2011), no. 3, 603–628. 

\medskip

[LP] J. Lebl and H. Peters, Polynomials constant on a hyperplane and CR maps of spheres,
 {\it Illinois J. Math.} 56 (2012), no. 1, 155-175 (2013). 

\medskip
[M] F. Meylan, Degree of a holomorphic map between unit balls from $C^2$ to $C^n$,  {\it Proc. Amer. Math. Soc.} 134 (2006), no. 4, 1023-1030.

\medskip

[P] S. Pinchuk, Proper holomorphic maps of strictly pseudoconvex domains, {\it Siberian Math. J.} 15 (1975), 644-649.

\medskip

[Q] Daniel G. Quillen, On the Representation of Hermitian Forms
as Sums of Squares, {\it Inventiones Math.} 5 (1968), 237-242.

\end{document}